\documentclass[oneside,english]{amsart}
\usepackage[T1]{fontenc}
\usepackage[latin9]{inputenc}
\usepackage{prettyref}
\usepackage{amstext}
\usepackage{amsthm}
\usepackage{amssymb}
\usepackage[numbers]{natbib}

\makeatletter
\numberwithin{equation}{section}
\numberwithin{figure}{section}
\theoremstyle{plain}
\newtheorem{thm}{\protect\theoremname}[section]
\theoremstyle{plain}
\newtheorem{lem}[thm]{\protect\lemmaname}
\theoremstyle{definition}
\newtheorem{example}[thm]{\protect\examplename}

\usepackage{mathtools}
\usepackage{braket}
\usepackage{prettyref}
\usepackage{tikz}
\usepackage{stmaryrd}
\usetikzlibrary{patterns}

\newrefformat{prop}{Proposition \ref{#1}}
\newrefformat{cor}{Corollary \ref{#1}}
\newrefformat{exa}{Example \ref{#1}}

\makeatother

\usepackage{babel}
\providecommand{\examplename}{Example}
\providecommand{\lemmaname}{Lemma}
\providecommand{\theoremname}{Theorem}

\begin{document}
\title{On the definition of neutrosophic logic}
\author{Takuma Imamura}
\address{Research Institute for Mathematical Sciences\\
Kyoto University\\
Kitashirakawa Oiwake-cho, Sakyo-ku, Kyoto 606-8502, Japan}
\email{timamura@kurims.kyoto-u.ac.jp}
\begin{abstract}
Smarandache (2003) introduced a new set-valued fuzzy logic called
(nonstandard) neutrosophic logic by using Robinson's nonstandard analysis.
However, its definition involved many errors including the illegal
use of nonstandard analysis. In this paper, we provide a rigorous
definition of neutrosophic logic. All the errors in the original definition
are addressed. We then point out some paradoxes of neutrosophic logic.
Finally we formulate neutrosophic logic with no use of nonstandard
analysis.
\end{abstract}

\maketitle

\section{Introduction}

\citet{Sma03} introduced a new set-valued fuzzy logic, called (nonstandard)
neutrosophic logic. This logic can be characterised as follows: each
proposition takes a value of the form $\left(T,I,F\right)$, where
$T$, $I$ and $F$ are subsets of the nonstandard unit interval $\left]\prescript{-}{}{0},1^{+}\right[$
(defined later) and represent all possible values of truthiness, indeterminacy
and falsity of the proposition, respectively. Unfortunately, the original
definition contains many errors including the illegal use of nonstandard
analysis.

In \prettyref{sec:Correction}, we first point out the errors that
appear in the original definition. We then provide a rigorous definition
of neutrosophic logic. All the errors in the original definition are
corrected. We also point out some problems concerning neutrosophic
logic, e.g., the paradox where complex propositions may have strange
truth values. In \prettyref{sec:Neutrosophic-logic-without}, we give
alternative definition of nonstandard neutrosophic logic without nonstandard
analysis.

Note that, in this paper, we do not discuss the theoretical significance
or the applications of neutrosophic logic. We only focus on mathematical
correctness of the target theory.

We assume the reader to have basic familiarity with elementary nonstandard
analysis (over the reals $\mathbb{R}$). We refer to \citet{Rob66}
for the foundations of nonstandard analysis.

\section{\label{sec:Correction}Correction of the definition}

We first review the original definition of neutrosophic logic given
in \citet{Sma03,Sma05} and keep track of the errors (and ambiguity).
We then provide a rigorous and corrected definition of neutrosophic
logic. Finally we point out some paradoxical (counter-intuitive) behaviours
of the neutrosophic logical connectives.

Before reviewing the definition, we recall the outline of the definition.
The usual fuzzy logic and its minor variations use the closed unit
interval $\left[0,1\right]_{\mathbb{R}}$ as the set of truth values.
Smarandache's (nonstandard) neutrosophic logic uses the nonstandard
unit interval $\left]\prescript{-}{}{0},1^{+}\right[$ of the hyperreals
$\mathbb{R}^{\ast}$ instead of the closed unit interval $\left[0,1\right]_{\mathbb{R}}$
of the reals $\mathbb{R}$. The nonstandard unit interval $\left]\prescript{-}{}{0},1^{+}\right[$
consists of three types of hyperreals: hyperreals between $0$ and
$1$, hyperreals infinitely close to but less than $0$, and hyperreals
infinitely close to but greater than $1$. Neutrosophic logic is then
characterised as the $\mathcal{P}\left(\left]\prescript{-}{}{0},1^{+}\right[\right)^{3}$\textendash valued
logic. Intuitively, for $\left(T,I,F\right)\in\mathcal{P}\left(\left]\prescript{-}{}{0},1^{+}\right[\right)^{3}$,
its components $T$, $I$ and $F$ represent the set of truthiness,
indeterminacy and falsity, respectively. The logical operations on
$\mathcal{P}\left(\left]\prescript{-}{}{0},1^{+}\right[\right)^{3}$
are defined according to this intuitive idea.

\subsection{Confused notations}

The notations $\prescript{-}{}{a}$ and $b^{+}$ are used as particular
hyperreal numbers.
\begin{quotation}
Let $\varepsilon>0$ be a such infinitesimal number. {[}...{]} Let's
consider the nonstandard finite numbers $1^{+}=1+\varepsilon$, where
``$1$'' is its standard part and ``$\varepsilon$'' its non-standard
part, and $\prescript{-}{}{0}=0-\varepsilon$, where ``$0$'' is its
standard part and ``$\varepsilon$'' its non-standard part. (\citep[p. 141]{Sma03};
\citep[p. 9]{Sma05})
\end{quotation}
At the same time, the notations $\prescript{-}{}{a}$ and $b^{+}$
are also used as particular sets of hyperreal numbers.
\begin{quote}
Actually, by ``$\prescript{-}{}{a}$'' one signifies a monad, i.e.,
a set of hyper-real numbers in non-standard analysis:
\[
\left(\prescript{-}{}{a}\right)=\set{a-x\in\mathbb{R}^{\ast}|x\text{ is infinitesimal}},
\]
and similarly ``$b^{+}$'' is a hyper monad:
\[
\left(b^{+}\right)=\set{b+x\in\mathbb{R}^{\ast}|x\text{ is infinitesimal}}.
\]
(\citep[p. 141]{Sma03}; \citep[p. 9]{Sma05})
\end{quote}
This confusion of notation can be found thereafter: in \citep[pp.10--11]{Sma05},
the notations $\prescript{-}{}{a}$ and $b^{+}$ are used as hypermonads;
on the other hand, in \citep[p.13]{Sma05}, $1^{+}$ is used as a
particular hyperreal number.

Note that the definitions of (one-sided) monads have minor errors.
The correct definitions are the following:
\begin{align*}
\left(\prescript{-}{}{a}\right) & :=\set{a-x\in\mathbb{R}^{\ast}|x\text{ is \emph{positive} infinitesimal}},\\
\left(b^{+}\right) & :=\set{b+x\in\mathbb{R}^{\ast}|x\text{ is \emph{positive} infinitesimal}}.
\end{align*}

\subsection{Ambiguous definition of the nonstandard unit interval}

\citet{Sma03} gives no explicit definition of the nonstandard unit
interval $\left]\prescript{-}{}{0},1^{+}\right[$ (or $\left\Vert \!\text{-}\prescript{-}{}{0},1^{+}\text{-}\!\right\Vert $
in \citep{Sma05}). He only says:
\begin{quotation}
Then, we call $\left]\prescript{-}{}{0},1^{+}\right[$ a non-standard
unit interval. Obviously, $0$ and $1$, and analogously non-standard
numbers infinitely small but less than $0$ or infinitely small but
greater than $1$, belong to the non-standard unit interval. (\citep[p. 141]{Sma03};
\citep[p. 9]{Sma05})
\end{quotation}
Here $\prescript{-}{}{0}$ and $1^{+}$ are particular real numbers
defined in the previous paragraph: $\prescript{-}{}{0}=0-\varepsilon$
and $1^{+}=1+\varepsilon$, where $\varepsilon$ is a fixed non-negative
infinitesimal. (Note that the phrase ``infinitely small but less than
$0$'' should be read as ``infinitely close to but less than $0$''.
Similarly, the phrase ``infinitely small but greater than $1$'' is
intended to mean ``infinitely close to but greater than $1$''.)

There are two possible definitions of the nonstandard unit interval:
\begin{enumerate}
\item $\left]\prescript{-}{}{0},1^{+}\right[=\set{x\in\mathbb{R}^{\ast}|\prescript{-}{}{0}<x<1^{+}}$
following the unusual notation of open interval;
\item $\left]\prescript{-}{}{0},1^{+}\right[=\set{x\in\mathbb{R}^{\ast}|0\substack{<\\
\approx
}
x\substack{<\\
\approx
}
1}$. (This is equal to the monad $\mu\left(\left[0,1\right]\right)$
of the closed unit interval.)
\end{enumerate}
In the first definition, it is false that hyperreal numbers infinitely
close to but less than $0$ or infinitely close to but greater than
$1$, belong to the nonstandard unit interval. The hyperreal $0-2\varepsilon$
is infinitely close to $0$ but not in $\left]\prescript{-}{}{0},1^{+}\right[$.
Similarly for $1+2\varepsilon$. Thus the first definition is not
compatible with the above-quoted sentence. The second definition is
better than the first one: when adopting the first one, the resulting
logic depends on the choice of the positive infinitesimal $\varepsilon$.
In our definition, we adopt the second one.

The cause of the confusion and the ambiguity can be found in the following
quote:
\begin{quotation}
We can consider $\left(\prescript{-}{}{a}\right)$ equals to the open
interval $\left(a-\varepsilon,a\right)$, where $\varepsilon$ is
a positive infinitesimal number. Thus:

$\left(\prescript{-}{}{a}\right)=\left(a-\varepsilon,a\right)$

$\left(b^{+}\right)=\left(b,b+\varepsilon\right)$

$\left(\prescript{-}{}{a}^{+}\right)=\left(a-\varepsilon_{1},a\right)\cup\left(a,a+\varepsilon_{2}\right)$,
where $\varepsilon,\varepsilon_{1},\varepsilon_{2}$ are positive
infinitesimal numbers. (\citep[p. 10]{Sma05})
\end{quotation}
Obviously it is wrong. Suppose, on the contrary, that $\left(\prescript{-}{}{a}\right)$
can be expressed in the form $\left(a-\varepsilon,a\right)$. Then
$a-\varepsilon$ does not belong to $\left(\prescript{-}{}{a}\right)$.
On the other hand, $a-\varepsilon$ is less than but infinitely close
to $a$, so $a-\varepsilon\in\left(\prescript{-}{}{a}\right)$, a
contradiction. This false belief well explains why Smarandache fell
into the confusion of notation and why he gave only an ambiguous description
of the nonstandard unit interval: if the monad could be described
like above, the two definitions of the unit interval would be equivalent.

\subsection{Misuse of nonstandard analysis}

Let us continue to read the definition.
\begin{quotation}
Let $T,I,F$ be standard or non-standard real subsets of $\left]\prescript{-}{}{0},1^{+}\right[$,

with $\sup T=t_{\sup}$, $\inf T=t_{\inf}$,

$\sup I=i_{\sup}$, $\inf I=i_{\inf}$,

$\sup F=f_{\sup}$, $\inf F=f_{\inf}$,

and $n_{\sup}=t_{\sup}+i_{\sup}+f_{\sup}$,

$n_{\inf}=t_{\inf}+i_{\inf}+f_{\inf}$.

The sets $T,I,F$ are not necessarily intervals, but may be any real
sub-unitary subsets: discrete or continuous; single-element, finite,
or (countably or uncountably) infinite; union or intersection of various
subsets; etc. (\citep[pp. 142--143]{Sma03}; \citep[p. 12]{Sma05})
\end{quotation}
Subsets of $\left]\prescript{-}{}{0},1^{+}\right[$ may have neither
infima nor suprema, because the transfer principle ensures the existences
of infima and suprema \emph{only for internal sets}. External sets
may lack suprema and/or infima. For instance, the monad $\mu\left(1/2\right)=\set{x\in\mathbb{R}^{\ast}|x\approx1/2}$
has neither the infimum nor the supremum. To see this, suppose, on
the contrary, that $\mu\left(1/2\right)$ has the infimum $L=\inf\mu\left(1/2\right)$.
Let $\varepsilon$ be any positive infinitesimal. Then there is an
$x\in\mu\left(1/2\right)$ such that $x\leq L+\varepsilon$. Since
every hyperreal infinitely close to $\mu\left(1/2\right)$ belongs
to $\mu\left(1/2\right)$, we have that $x-2\varepsilon\in\mu\left(1/2\right)$.
Hence $L\leq x-2\varepsilon\leq L-\varepsilon$, a contradiction.
The same applies to the supremum.

There are two workarounds:
\begin{enumerate}
\item inserting the sentence ``assume that $T,I,F$ are internal'' or ``assume
that $T,I,F$ have infima and suprema''; or
\item giving up the use of infima/suprema.
\end{enumerate}
When we adopt the first one, we cannot consider propositions of external
values. For example, the whole interval $\left]\prescript{-}{}{0},1^{+}\right[$
is external, so we cannot consider \emph{completely ambiguous} propositions,
none of whose truthiness, indeterminacy and falsity are (even roughly)
determined. In fact, none of infima and suprema are necessary to formulate
neutrosophic logic. In our definition, we adopt the second prescription.
So we can consider propositions with external values such as $\left]\prescript{-}{}{0},1^{+}\right[$.

\subsection{Rigorous definition of neutrosophic logic}

Now let us correct the definition of neutrosophic logic. We first
define the nonstandard unit interval as follows:
\[
\left]\prescript{-}{}{0},1^{+}\right[=\set{x\in\mathbb{R}^{\ast}|0\substack{<\\
\approx
}
x\substack{<\\
\approx
}
1}.
\]
Let $\mathbb{V}$ be the power set of $\left]\prescript{-}{}{0},1^{+}\right[$,
the collection of all subsets of $\left]\prescript{-}{}{0},1^{+}\right[$.
We define the binary operations $\owedge$, $\ovee$ and $\obslash$
on $\mathbb{V}$ as follows:
\begin{align*}
A\owedge B & =\set{ab|a\in A,b\in B},\\
A\ovee B & =\set{a+b-ab|a\in A,b\in B},\\
A\obslash B & =\set{c-a+ab|a\in A,b\in B,c\in1^{+}}.
\end{align*}
We need to verify the closure property of $\mathbb{V}$.
\begin{lem}
$\mathbb{V}$ is closed under $\owedge$, $\ovee$ and $\obslash$.
\end{lem}

\begin{proof}
Consider the following standard operations on $\mathbb{R}$:
\begin{align*}
f\left(a,b\right) & =ab,\\
g\left(a,b\right) & =a+b-ab,\\
h\left(c,a,b\right) & =c-a+ab.
\end{align*}
Notice that they are continuous everywhere. The standard unit interval
$\left[0,1\right]$ is closed under $f,g$ and $h\left(1,\cdot,\cdot\right)$.
Here we only prove the case of $h\left(1,\cdot,\cdot\right)$. Since
$h\left(1,a,b\right)$ is monotonically decreasing for $a$, we have
that $\min_{\left(a,b\right)\in\left[0,1\right]^{2}}h\left(1,a,b\right)=\min_{b\in\left[0,1\right]}h\left(1,1,b\right)=0$.
Similarly, since $h\left(1,a,b\right)$ is monotonically increasing
for $b$, we have that $\max_{\left(a,b\right)\in\left[0,1\right]^{2}}h\left(1,a,b\right)=\max_{a\in\left[0,1\right]}h\left(1,a,1\right)=1$.
Thus $h\left(\set{1}\times\left[0,1\right]\times\left[0,1\right]\right)\subseteq\left[0,1\right]$.

Now, let $a,b\in\left]\prescript{-}{}{0},1^{+}\right[$ and $c\in1^{+}$.
Choose $a',b'\in\left[0,1\right]$ infinitely close to $a,b$, respectively.
Of course, $c$ is infinitely close to $1$. By the nonstandard characterisation
of continuity (see Theorem 4.2.7 of \citep{Rob66}), $f^{\ast}\left(a,b\right),g^{\ast}\left(a,b\right),h^{\ast}\left(c,a,b\right)$
are infinitely close to $f\left(a',b'\right),g\left(a',b'\right),h\left(1,a',b'\right)\in\left[0,1\right]$,
respectively. Hence $f^{\ast}\left(a,b\right),g^{\ast}\left(a,b\right),h^{\ast}\left(c,a,b\right)\in\left]\prescript{-}{}{0},1^{+}\right[$.

Let $A,B\in\mathbb{V}$. Then $A\owedge B=f^{\ast}\left(A\times B\right)$,
$A\ovee B=g^{\ast}\left(A\times B\right)$ and $A\obslash B=h^{\ast}\left(1^{+}\times A\times B\right)$
belong to $\mathbb{V}$.
\end{proof}
According to the original definition (\citep[p. 143]{Sma03}), neutrosophic
logic is the $\mathbb{V}^{3}$-valued (extensional) logic. Each proposition
takes a value of the form $\left(T,I,F\right)\in\mathbb{V}^{3}$,
where $T$ represents possible values of truthiness, $I$ indeterminacy,
and $F$ falsity. The logical connectives $\wedge,\vee,\to$ are defined
as follows:
\begin{align*}
\left(T_{1},I_{1},F_{1}\right)\wedge\left(T_{2},I_{2},F_{2}\right) & =\left(T_{1}\owedge T_{2},I_{1}\owedge I_{2},F_{1}\owedge F_{2}\right),\\
\left(T_{1},I_{1},F_{1}\right)\vee\left(T_{2},I_{2},F_{2}\right) & =\left(T_{1}\ovee T_{2},I_{1}\ovee I_{2},F_{1}\ovee F_{2}\right),\\
\left(T_{1},I_{1},F_{1}\right)\to\left(T_{2},I_{2},F_{2}\right) & =\left(T_{1}\obslash T_{2},I_{1}\obslash I_{2},F_{1}\obslash F_{2}\right).
\end{align*}

Note that our definitions of $\ovee$ and $\obslash$ are different
from the original ones. Smarandache uses the following operations
instead of $\ovee$ and $\obslash$ (\citep[p. 145]{Sma03}):
\begin{align*}
A\ovee'B & =\left(A\oplus B\right)\ominus\left(A\odot B\right),\\
A\obslash'B & =1^{+}\ominus A\oplus\left(A\odot B\right),
\end{align*}
where $\ominus,\odot,\oplus$ are the elementwise subtraction, multiplication
and addition of sets. There are at least two reasons why the original
definition is not good. The first is pre-mathematical. When calculating,
for example, $A\ovee'B=\left(A\oplus B\right)\ominus\left(A\odot B\right)$,
the second and the third occurrences of $A$ can take different values,
despite that they represent the same proposition. The same applies
to $B$ and the calculation of $A\obslash'B$ obviously. The second
is mathematical. $\mathbb{V}$ is not closed under those operations:
\begin{align*}
2 & =\text{\ensuremath{\underbar{1}}}+1-\text{\ensuremath{\underbar{0}}}\cdot1\in\set{0,1}\ovee'\set{1},\\
2+\varepsilon & =1+\varepsilon-\text{\ensuremath{\underbar{0}}}+\text{\ensuremath{\underbar{1}}}\cdot1\in\set{0,1}\obslash'\set{1},
\end{align*}
where $\varepsilon$ is positive infinitesimal. Because of this, the
following ad-hoc workaround is needed:
\begin{quotation}
{[}...{]} if, after calculations, one obtains number $<0$ or $>1$,
one replaces them by $\prescript{-}{}{0}$ or $1^{+}$, respectively.
(\citep[p. 145]{Sma03})
\end{quotation}

\subsection{Paradoxical phenomena}

Consider the $\mathbb{V}$\textendash valued logic, where each proposition
takes a \emph{truth} value $T\in\mathbb{V}$. Each neutrosophic logical
connectives was defined componentwise. In other words, neutrosophic
logic is the $3$\textendash fold product of the $\mathbb{V}$\textendash valued
logic. Hence neutrosophic logic cannot be differentiated from the
$\mathbb{V}$\textendash valued logic by equational properties (see
\citet[Lemma 11.3 of Chapter II]{BS81}).

This causes some paradoxical phenomena. Let $A$ be a (true) proposition
with value $\left(\set{1},\set{0},\set{0}\right)$ and let $B$ be
a (false) proposition with value $\left(\set{0},\set{0},\set{1}\right)$.
Usually we expect that the falsity of the conjunction $A\wedge B$
is $\set{1}$. However, its actual falsity is $\set{0}$. We expect
that the indeterminacy of the negation $\neg A$ is $\set{0}$. However,
its actual indeterminacy is $1^{+}$ (see \citet[p. 145]{Sma03} for
the definition of the negation).

This problem has been addressed in \citet{Riv08}.

\section{\label{sec:Neutrosophic-logic-without}Neutrosophic logic without
nonstandard analysis}

\subsection{Nonarchimedean fields}

Let $\mathbb{K}$ be an ordered field. It is well-known that the ordered
semiring $\mathbb{N}$ of natural numbers can be canonically embedded
into $\mathbb{K}$ by sending $n\mapsto\underset{n}{\underbrace{1_{\mathbb{K}}+\cdots+1_{\mathbb{K}}}}$.
This embedding can be uniquely extended to the ordered ring $\mathbb{Z}$
of integers and to the ordered field $\mathbb{Q}$ of rational numbers.
Thus we may assume without loss of generality that $\mathbb{Q}\subseteq\mathbb{K}$.
An element $x\in\mathbb{K}$ is said to be \emph{infinitesimal} (relative
to $\mathbb{Q}$) if for any positive $q\in\mathbb{Q}$, $\left|x\right|\leq q$.
For instance, the unit of the addition $0_{\mathbb{K}}$ is trivially
infinitesimal. The ordered field $\mathbb{K}$ is called \emph{nonarchimedean}
if it has nonzero infinitesimals.
\begin{example}
Every ordered subfield of the real field $\mathbb{R}$ is archimedean.
Conversely, every archimedean ordered field can be (uniquely) embedded
into $\mathbb{R}$ (see \citet[Theorem 10.21]{Bly05}).
\end{example}

\begin{example}
The hyperreal field $\mathbb{R}^{\ast}$ is a nonarchimedean ordered
field. Generally, every proper extension $\mathbb{K}$ of $\mathbb{R}$
is nonarchimedean: Let $x\in\mathbb{K}\setminus\mathbb{R}$. There
are two cases. Case I: $x$ is \emph{infinite} (i.e. its absolute
value $\left|x\right|$ is an upper bound of $\mathbb{R}$). Then
its reciprocal $1/x$ is nonzero infinitesimal. Case II: $x$ is finite.
Then the set $\set{y\in\mathbb{R}|y<x}$ is nonempty and bounded in
$\mathbb{R}$. So it has the supremum $x^{\circ}$. The difference
$x-x^{\circ}$ is nonzero infinitesimal. Hence $\mathbb{K}$ is nonarchimedean.
\end{example}

\begin{example}[{cf. \citet[pp. 15--16]{GO03}}]
Let $\mathbb{K}$ be an ordered field. Let $\mathbb{K}\left(X\right)$
be the field of rational functions over $\mathbb{K}$. Define an ordering
on $\mathbb{K}\left(X\right)$ by giving a positive cone:
\[
0\leq\frac{f\left(X\right)}{g\left(X\right)}\iff0\leq\frac{\text{the leading coefficients of }f\left(X\right)}{\text{the leading coefficients of }g\left(X\right)}.
\]
Then $\mathbb{K}\left(X\right)$ forms an ordered field having nonzero
infinitesimals (relative to not only $\mathbb{Q}$ but also $\mathbb{K}$)
such as $1/X$ and $1/X^{2}$. Hence $\mathbb{K}\left(X\right)$ is
nonarchimedean.
\end{example}

\subsection{Alternative definition of neutrosophic logic}

Comparing with other nonarchimedean fields, one of the essential features
of the hyperreal field $\mathbb{R}^{\ast}$ is the transfer principle,
which states that $\mathbb{R}^{\ast}$ has the same first-order properties
as $\mathbb{R}$. On the other hand, the formulation of neutrosophic
logic does not depend on the transfer principle. The use of nonstandard
analysis is not essential for this logic, and can be eliminated from
its definition.

Fix a nonarchimedean ordered field $\mathbb{K}$. For $x,y\in\mathbb{K}$,
$x$ and $y$ are said to be \emph{infinitely close} (denoted by $a\approx b$)
if $a-b$ is infinitesimal. We say that $x$ is \emph{roughly smaller
than} $y$ (and write $x\substack{<\\
\approx
}
y$) if $x<y$ or $x\approx y$. For $a,b\in\mathbb{K}$ the set $\left]\prescript{-}{}{a},b^{+}\right[_{\mathbb{K}}$
is defined as follows:
\[
\left]\prescript{-}{}{a},b^{+}\right[_{\mathbb{K}}=\set{x\in\mathbb{K}|a\substack{<\\
\approx
}
x\substack{<\\
\approx
}
b}.
\]
Let $\mathbb{V_{K}}$ be the power set of $\left]\prescript{-}{}{0},1^{+}\right[_{\mathbb{K}}$.
The $\mathbb{K}$\textendash valued neutrosophic logic is defined
as the $\mathbb{V_{K}}^{3}$\textendash valued logic. The rest of
the definition is completely the same as the case $\mathbb{K}=\mathbb{R}^{\ast}$.

\bibliographystyle{IEEEtranSN}
\bibliography{bib}

\end{document}